\documentclass[11pt]{amsart}
\usepackage{amssymb}
\voffset- .3cm
\numberwithin{equation}{section}

\newcommand\A{\mathcal{A}}

\newcommand{\N}{\mathbb{N}}
\newcommand{\R}{\mathbb{R}}
\newcommand{\C}{\mathbb{C}}
\newcommand{\Z}{\mathbb{Z}}

\usepackage{amsmath,amssymb}
\usepackage{picins,slashbox,dbnsymb,
  graphicx,verbatim,longtable,
  epic,eepic,stmaryrd,rotating}
\usepackage[all]{xy}

\theoremstyle{plain}

\theoremstyle{definition}

\theoremstyle{remark}

\newlength{\standardunitlength}
\catcode`\@=11
\long\def\@makecaption#1#2{%
    \vskip 10pt
    \setbox\@tempboxa\hbox{
      \small{#1: }\ignorespaces #2}%
    \ifdim \wd\@tempboxa >\captionwidth {%
        \rightskip=\@captionmargin\leftskip=\@captionmargin
        \unhbox\@tempboxa\par}%
      \else
        \hbox to\hsize{\hfil\box\@tempboxa\hfil}%
    \fi}

\newdimen\@captionmargin\@captionmargin=2\parindent
\newdimen\captionwidth\captionwidth=\hsize
\catcode`\@=12

\newlength{\globalparindent}
\begin{document}

\hfill\hfill To my son Philippe for his unbounded
energy and optimism.

\vskip 5cm

\centerline {\bf On the classification of Floer-type theories.}

\vskip 1cm

\centerline {\bf Nadya Shirokova.}

\vskip 2cm

\centerline {\bf Abstract}

\vskip 1cm

 In this paper we outline a program for the classification of  Floer-type theories, (or defining 
 invariants of finite type for families). We  consider Khovanov complexes as a local
system on the space of knots introduced by V. Vassiliev and construct the wall-crossing morphism. We extend this system
to the singular locus by the cone of this morphism and introduce the definition of the local system of finite type. 
This program can be further generalized to the manifolds of dimension 3 and 4 [S2], [S3].

\vskip 1cm

\newpage

\vskip 3cm

\centerline {\bf Contents}

\vskip.3cm
{\bf 1. Introduction.}
\vskip .5cm

{ \bf 2. Vassiliev's and Hatcher's theories.}
\vskip.4cm

2.1. The space of knots, coorientation. 
\vskip.2cm
2.2. Vassiliev derivative.

\vskip.2cm
2.3. The topology of the chambers of the space of knots.

\vskip .5cm
{ \bf 3. Khovanov homology.}
\vskip.5cm
3.1. Jones polynomial as Euler characterictics.  Skein relation.
\vskip.2cm
3.2. Reidemeister and Jacobsson moves.
\vskip.2cm
3.3. Wall-crossing morphisms.
\vskip .2cm
3.4. The local system of Khovanov complexes on the space of knots.

\vskip .5cm
{\bf 4. Main definition, invariants of finite type for families.}
\vskip.4cm

 4.1. Some homological algebra.
\vskip.2cm
4.2. Space of knots and the classifying space of the category.
\vskip.2cm
4.3. Vassiliev derivative as a cone of the wall-crossing morphism.
\vskip.2cm
4.4. The definition of a theory of finite type.

\vskip .5cm
{\bf 5. Theories of finite type. Further directions.}
\vskip.4cm
5.1. Examples of combinatorially defined theories.
\vskip.2cm
5.2.  Generalizations to dimension 3 and 4.
\vskip.2cm
5.3. Further directions.
\vskip.5cm

{\bf 6. Bibliography}.
\vskip 1.5cm
\vskip .5cm

\newpage
\centerline {\bf 1. Introduction.}

\vskip 1cm

    Lately there has been a lot of interest in various categorifications of classical scalar invariants,
  i.e. homological theories, Euler characteristics of which are scalar invariants. Such examples
  include the original instanton Floer homology, Euler characteristic of which, as it was proved
  by C.Taubes [T], is Casson's invariant. Ozsvath-Szabo [OS] 3-manifold theory
  categorifies Turaev's torsion, the Euler characteristic of their knot homologies [OS] is the Alexander   polynomial. The theory of M. Khovanov categorifies the Jones polynomial  [Kh] and Khovanov-Rozhansky theory categorifies the $sl(n)$ invariants [KR] .
    \vskip .2cm
  
   The theory that we are constructing will bring together theories of V. Vassiliev, A. Hatcher and
   M. Khovanov, and while describing their results we will specify which parts of their constructions
   will be important to us.
  
    The resulting theory can be considered as a "categorification of Vassiliev
   theory" or a classification of categorifications of knot invariants. We introduce the definition of a 
   theory of finite type n and show that Khovanov homology theory in a categorical sense
    decomposes into a "Taylor series" of theories of finite type.

     The Khovanov functor is just the first example of a theory satisfying our axioms
   and we believe, that all theories mentioned above will fit into our template.
    \vskip .2cm
     Our  main strategy is to consider a knot homology theory as a local system, or a constructible
 sheaf on the space of all objects (knots, including singular ones), extend this local system to the singular locus and introduce the analogue of the "Vassiliev derivative" for categorifications.

 By studying spaces of embedded manifolds we implicitly study their diffeomorphism groups and
invariants of finite type. In his seminal paper [V] Vassiliev introduced  finite type invariants by considering the space
 of all immersions of $S^1$ into $R^3$ and relating  the topology of
 the singular locus to the topology of its complement via Alexander duality. He resolved and cooriented the discriminant of the space and  introduced a
 spectral sequence with a filtration, which suggested the simple geometrical  and combinatorial definition of an invariant of finite type, which
was later interpreted by Birman and Lin as a "Vassiliev derivative" and led to
the following skein relation.

If $\lambda$ be an arbitrary invariant of oriented knots in oriented space
with values in some Abelian group $A$. Extend $\lambda$ to be an invariant of
$1$-singular knots (knots that may have a single singularity that locally
looks like a double point $\doublepoint$), using the formula

$$  \lambda(\doublepoint)=\lambda(\overcrossing)-\lambda(\undercrossing)$$

Further extend $\lambda$ to the set  of $n$-singular knots (knots
with $n$ double points) by repeatedly using the skein relation.
    \vskip .2cm
{\bf Definition} We say that $\lambda$ is of type $n$ if its extension
to $(n+1)$-singular knots vanishes
identically. We say that $\lambda$ is of finite type if it is of type $n$ for
some $n$.
    \vskip .2cm

 Given the above formula, the definition of an invariant of finite type n becomes similar to that  of a polynomial:  its (n+1)st
 Vassiliev derivative is zero.
 
 It was shown that all known invariants are either of finite type, or are infinite linear combinations of those,
 e.g. in  [BN1] it was shown that the nth coeffitient of the Conway polynomial is a Vassiliev invariant of order $\leq n$.

 In this paper we are working with  Khovanov homology, which will be
 our main example, however the latest progress in finding the combinatorial formula  for the differential 
 of the  Ozsvath-Szabo knot complex [MOS], makes us hopeful that more and more examples will be coming.

 For the construction of the local system it is important to understand the topological type of the base. 
 The topology of the connected components of the complement to the discriminant in the space of knots, called chambers, was studied by A. Hatcher and R.Budney  [H], [B].
They introduced  simple homotopical models for such spaces. Recall that the local system is well-defined
on a homotopy model of the base, so Hatcher's model is exactly what is needed to construct the
local system of Khovanov complexes.

Throughout the paper the following observation is the main guideline for our constructions:
{\sl local systems of the classifying space of the category are functors from this category to the triangulated category of complexes}.

 It would be very interesting to understand the relation between the
Vassiliev space of knots is the classifying space of the category, whose objects  are knots and whose morphisms are knot cobordisms.

 Our construction provides a {\sl Khovanov functor} from the category of knots into the triangulated category of complexes.

This  allows us to translate all topological properties of the space of knots and Khovanov
local system on it into the language of homological algebra and then use the methods of triangulated categories and 
homological algebra to assign algebraic objects to topological ones (singular knots and links).

  \vskip .3cm
    Recall that in  his paper [Kh]  M.Khovanov categorified the Jones polynomial, i.e. he found a homology theory, the
  Euler characteristics of which equals the Jones polynomial.   He starts with a diagram of the knot and
  constructs a bigraded complex, associated to this diagram, using two resolutions of the knot crossing:
  
 \vskip .4cm
 \begin{center}
 \includegraphics{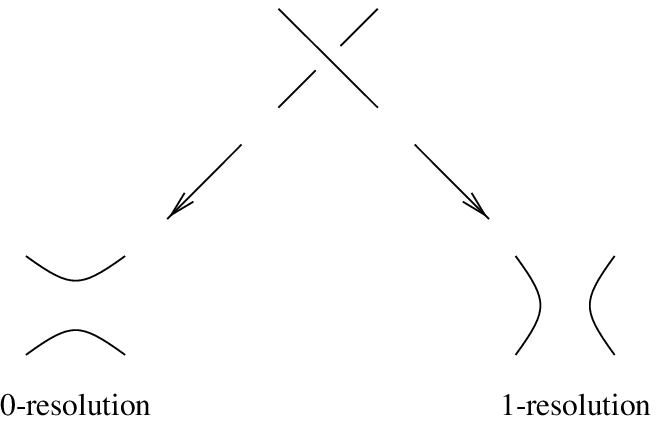}
 \end{center}
 \vskip .4cm
  
  The Khovanov complex then becomes the sum of the tensor products of the vector space V, where the
  homological degree is given by the number of 1's in the complete resolution of the knot.
   The
 local system of Khovanov homologies on the Vassiliev's space of knots can be considered as invariants of families
 of knots.
 
  The discriminant of  Vassiliev's space corresponds to knots with transversal self-intersection, i.e. moving
  from one chamber to another we change overcrossing to undercrossing by passing through a knot
  with a single double point. We study how the Khovanov complex changes under such modification and find the corresponding morphism.

  After defining a wall-crossing morphism we can extend the  invariant
to the singular locus by the cone of a morphism which is our "categorification of the Vassiliev derivative". Then we introduce the definition of a local system of finite type: the local system is of finite type n if for any selfintersection of the discriminant of codimension n, its n'th cone is an acyclic complex.

  \vskip .3cm

    The categorification of the Vassiliev derivative allows us to define the filtration on the Floer - type theories for manifolds of any dimension.

  In [S4] we prove the first finiteness result:

\vskip.3cm

{\bf Theorem [S4]}. Restricted to the subcategory of knots with at most $n$ crossing, Khovanov local system is of finite type $n$, $n\geq 3$ and of type zero $n=0,1,2.$
\vskip.3cm

   This definition can be generalized to the categorifications of the invariants of manifolds of any dimension:  we construct spaces of 3 and 4-manifolds by a version of a Pontryagin-Thom construction,
   consider homological invariants of 3 and 4-manifolds as local systems on these spaces and extend
   them to the discriminant.

  In subsequent  papers   our main example will be the Heegaard Floer homology [OS], the Euler characteristic of which is Turaev's torsion.   We show that local systems of such homological theories on the space of 3 - manifolds [S1] will carry information about invariants of finite type for families and information about the diffeomorphism group.
   We also have a construction [S2] for the  refined Seiberg-Witten invariants on the
space of parallelizable 4-manifolds.
  \vskip.3cm
{\bf Acknowledgements}. My deepest thanks go  to Yasha Eliashberg  for many valuable discussions,
for inspiration and
for his constant encouragement and  support. I want to thank Maxim Kontsevich  who suggested that I work on this project, for his  attention to my work during my visit to the IHES and for many  important suggestions.

 I want to thank graduate students Eric Schoenfeld and Isidora Milin for reading the paper and making useful comments. 
 
 This paper was written during my visits to the IAS, IHES, MPIM and Stanford and I am grateful to these institutions for their exceptional hospitality. This work was partially supported by the NSF grant DMS9729992.
  
\vskip .5cm
\newpage

\centerline{ \bf 2. Vassiliev theory, invariants of finite type.}
\vskip 1cm
{ \sl 2.1. The space of knots, coorientation}
\vskip.5cm
 Vassiliev considered the space of all maps  $E = f: S^1 \rightarrow R^3$.
 This space is  a space of functions, so it is an infinite-dimensional Euclidean
 space. It is linear, contractible, and consists of singular (D) and nonsingular(E - D) knots. The 
 discriminant D forms a singular hypersurface in E and subdivides into chambers, corresponding
 to different isotopy types of knots. To move from one chamber to another one has to change one overcrossing to undercrossing, passing through a singular knot with one double point.
 \vskip.2cm
 
 The discriminant of the space of knots is a real hypersurface, stratified by the number of the double points, which subdivides the infinite-dimensional 
 space into {\sl chambers}, corresponding to different isotopy types of knots.
 
 Vassiliev resolved and cooriented the discriminant, so we can assume that all points of selfintersection
 are transversal, with $2^n$ chambers adjacent to a point of
 selfintersection of the discriminant of codimension $n$.

 To study the topology of the complement to the discriminant, Vassiliev wrote a spectral sequence,
 calculating the homology of the discriminant and then related it to the homology of its complement
 via Alexander duality.  His
 spectral sequence had a filtration, which suggested the simple geometrical  and combinatorial definition of an invariant of finite type: an invariant is of type $n$ if for any selfintersection of the discriminant of codimension (n+1) its alternated sum over the $2^{n+1}$ chambers adjacent to a point of
 selfintersection is zero.

\vskip.4cm

  For our constructions it will be very important to have a coorientation of the discriminant, which
  was introduced by Vassiliev.
  
  \vskip.2cm
{\bf Definition}. A hypersurface in a real manifold is said to be {\sl coorientable} if it has a non-zero section of its normal bundle, i.e. if there exists a continuous vector field which is not tangent to the hypersurface at any point and doesn't vanish anywhere.  
\vskip .2cm

So there are two sides of the hypersurface : one where this vector field is pointing to and the other is where it is pointing from. And there are two choices of such vector field. The {\sl coorientation} of a coorientable hypersurface is the choice of one of two possibilities.

 For example, Mobius band in $R^3$ is not coorientable.
 \vskip.2cm
 Vassiliev shows [V] that the discriminant of the space of knots has a coorientation, the conistent choice of normal directions.
 
 Recall that the nonsingular point  $ \psi \in D$ of the discriminant is a map  $S^1 \rightarrow R^3$, gluing together
 2 distinct points $t_1, t_2$ of  $S^1$, s.t. derivatives of the map $\psi$ at those points are transversal.
 \vskip.4cm
 {\bf Coorientation of the discriminant}. Fix the orientation of   $R^3$ and choose positively oriented local coordinates near the point
 $\psi(t_1) = \psi(t_2)$. For any point $\psi_1 \in D$ close to $\psi$ define the number $r(\psi_1)$
 as the determinant:
 $$(\frac {\partial \psi_1} {\partial \tau} \arrowvert t_1,\frac {\partial \psi_1} {\partial \tau} \arrowvert t_2, \psi_1(t_1) - \psi_1(t_2))$$
 
 with respect to these coordinates. This determinant depends only of the pair of points $t_1, t_2$,
 not on their order. A vector in the space of functions at the point   $ \psi \in D$, which is transversal
 to the discriminant, is said to be positive, if the derivative of the function r along this vector is positive and negative,
 if this derivative is negative.
 
 This rule gives the coorientation of the hypersurface $D$ at all its nonsingular points and also of any
 nonsingular locally irreducible component of D at the points of selfintersection of D.
 \vskip .2cm
 
 The consistent choice of the normal directions of the walls of the discriminant will give the "directions"
  of the
 cobordisms (which are embedded into $E \times I$) between knots of the space E.
\vskip .2cm

  \vskip.2cm
{\bf Note}. It is interesting to compare this construction with the result of E.Ghys [Gh], who introduced a metric on the space of knots and 3-manifolds.)
\vskip .5cm

{\sl 2.3. The topology of the chambers of the space of knots.}
\vskip .5cm

 The study of the topology of the chambers of the space of knots was started by A. Hatcher [H], who found
 a simple homotopy models for these spaces.

\vskip .2cm
 
  The main result is based on an earlier theorem regarding the topology of the classifying space of
 diffeomorphisms of an irreducible 3-manifold with nonempty  boundary.
 
 In the following theorem A. Hatcher and D. McCullough answered the question posed by M. 
Kontsevich [K], regarding the finiteness of the homotopy type of the classifying space of the 
group of diffeomorphisms [HaM]: 
\vskip .2cm

{\bf Theorem [HaM].} Let M be an irreducible compact connected orientable 3-manifold with 
nonempty boundary. Then $BDif f (M, rel\partial)$ has the homotopy type of a finite aspherical CW- 
complex.
    \vskip .2cm
The proof of this theorem uses the JSJ-decomposition of a 3-manifold.

When applied to knot complements, 
the JSJ-decomposition defines a fundamental
class of links in $S^3$, the "knot generating links" (KGL).
A KGL is any $(n+1)$-component link 
$L=(L_0,L_1,\cdots,L_n)$ whose complement is either
Seifert fibred or atoroidal, such that the $n$-component
sub-link $(L_1,L_2,\cdots,L_n)$ is the unlink. 
 If the complement
of a knot $f$ contains an incompressible torus, then 
$f$ can be represented as a `spliced knot' $f=J \Box L$
in  unique way, where $L$ is an $(n+1)$-component KGL,
 and 
$J=(J_1,\cdots,J_n)$ is an $n$-tuple
of non-trivial long knots.

 The spliced knot  $J \Box L$ is obtained from $L_0$ by a
generalized satellite construction.
For any knot there is a representation of a knot as an iterated
splice knot of atoroidal and hyperbolic KGLs.  The order of
splicing determines the "companionship tree"
of $f$, $G_f$, and is a complete isotopy invariant of long knots.

Given a knot $f \in K$, denote the path-component of $K$ containing $f$
by $K_f$. 
 The topology of the chambers $K_f$ was further studied by R. Budney
 The main result of his paper [Bu]
is the computation of the homotopy type of
$K_f$ if $f$ is a hyperbolically-spliced knot ie: $f=J \Box L$ where
$L$ is a hyperbolic KGL. 
\vskip .2cm
The combined results can
 be summarized in the following theorem:
 
 \vskip .2cm
 \newpage
 {\bf Theorem [Bu, H]}.  
 \vskip 1cm

If $f=J \Box L$ where $L$ is an $(n+1)$-component hyperbolic KGL, then 
$$K_f \backsimeq S^1 \times \left( SO_2 \times_{A_f} \prod_{i=1}^n K_{J_i} \right)$$
$A_f$ is the maximal subgroup of $B_L$ such that induced action of 
$A_f$ on $K^n$ preserves $\prod_{i=1}^n K_{L_i}$.
The restriction map $A_f \to Diff(S^3,L_0) \to Diff(L_0)$
is faithful, giving an embedding $A_f \to SO_2$, and this is the
action of $A_f$ on $SO_2$.

This result completes the computation of the homotopy-type of $K$ since
we have the prior results:
 
\begin{enumerate}
 \item [H1] If $f$ is the unknot, then $K_f$ is contractible.
 \item [H2] If $f$ is a torus knot, then $K_f \simeq S^1$.
 \item [H3] If $f$ is a hyperbolic knot, then 
  $K_f \backsimeq S^1 \times S^1$
 \item [H4] If a knot $f$ is a cabling of a 
 knot $g$ then $K_f \backsimeq S^1 \times K_g$. 
 \item [B5] If the knot $f$ is a connected sum of
 $n \geq 2$ prime knots $f_1, f_2, \cdots, f_n$ then
 $K_f \backsimeq \left( ({\mathcal C}_2(n) \times \prod_{i=1}^n K_{f_i}\right)/\Sigma_f$.
 Here $\Sigma_f \subset S_n$ is a Young subgroup of $S_n$, acting on
 ${\mathcal C}_2(n)$ by permutation of the labellings of the cubes, and similarly by permuting
 the factors of the product $\prod_{i=1}^n K_{f_i}$.  The definition of
 $\Sigma_f \subset S_n$ is that it is the subgroup of $S_n$ that preserves a
 partition of $\{1,2,\cdots,n\}$, the partition being given by the equivalence
 relation $i \sim j \Longleftrightarrow K_{f_i} =K_{f_j}$.
 \item [B6] If a knot has a non-trivial companionship tree, then
  it is either a cable, in which case H4 applies, 
  a connect-sum, in which case B5 applies or is hyperbolically
  spliced.
  If
  a knot has a trivial companionship tree, it is either the unknot,
  in which case H1 applies, or a torus knot in which case
  H2 applies, or a hyperbolic knot, in which case H3 applies.
  Moreover, every time one applies one of the above theorems,
  one reduces the problem of computing the homotopy-type of 
  $K_f$ to computing the homotopy-type of knot spaces for knots 
  with shorter companionship trees,
  thus the process terminates after finitely-many iterations.
\end{enumerate}

\vskip .5cm
 For constructing a local system we need only the homotopy type of the chamber. The theorem of Hatcher and Budney provides us with a complete classification of homotopy types of chambers, corresponding to all possible knot types.

\newpage
\centerline{ \bf 3. Khovanov's categorification of Jones polynomial.}

  \vskip 1cm

{ \sl 3.1. Jones polynomial as Euler characterictics.  Skein relation.}

\vskip .5cm
   In his paper [Kh]  M. Khovanov constructs a homology theory, with Euler characteristics 
   equal to the
   Jones polynomial.

 \vskip .3cm
 He associated to any diagram $D$ of an oriented 
link with n crossing points a chain complex $CKh(D)$ of abelian groups of homological length $(n+1)$,
and proved that for 
any two diagrams of the same link the corresponding 
complexes are chain homotopy equivalent. Hence, the 
homology groups $Kh(D)$ are link invariants up to isomorphism.

  His construction is as follows: given any double point of the link projection $D$, he allows two smoothings:

 \vskip .4cm
 \begin{center}
 \includegraphics{ch4p5.eps}
 \vskip .4cm
 \end{center}

 If the the diagram has n double points, there are $2^n$ possible resolutions. The result of each 
 complete smoothing is the set of circles in the plane, labled by $n$-tuples of 1's and 0's:

 $$CKh( \underbrace{\bigcirc,...,\bigcirc}_{n times}) = V^{\otimes n}$$

   The cobordisms between links, i.e.,
   surfaces embedded in $\R^3\times [0,1],$ 
   should provide maps between the associated 
  groups.   A surface embedded in the 4-space can be visualized as a sequence 
  of plane projections of its 3-dimensional sections (see [CS]). 
  Given such a presentation $J$ of a compact oriented surface $S$ 
  properly embedded in $\R^3\times [0,1]$ with the boundary of $S$ 
  being the union of two links $L_0\subset \R^3\times \{ 0\} $ and 
  $L_1 \subset \R^3\times \{ 1\},$  , Khovanov
   associates to $J$ a map 
  of cohomology groups 
  $$
  \theta_J: Kh^{i,j}(D_0)\rightarrow Kh^{i, j + \chi(S)}(D_1), \hspace{0.4in} 
  i,j\in \Z
 $$
  
   The differential of the Khovanov complex is defined using two
linear maps $m:V\otimes V\to V$ and $\Delta:V\to V\otimes V$ 
given by formulas :

\vskip 1cm
  $$  \big(V\otimes V\overset{m}{\rightarrow}V\big)
    \quad m:\begin{cases}
      v_+\otimes v_-\mapsto v_- &
      v_+\otimes v_+\mapsto v_+ \\
      v_-\otimes v_+\mapsto v_- &
      v_-\otimes v_-\mapsto 0
    \end{cases}$$
$$    \big(V\overset{\Delta}{\rightarrow}V\otimes V\big)
    \quad \Delta:\begin{cases}
      v_+ \mapsto v_+\otimes v_- + v_-\otimes v_+ &\\
      v_- \mapsto v_-\otimes v_- &
    \end{cases}$$

\vskip 1cm
 The differential in Khovanov complex can be informally described as "all the ways of changing 0-crossing to 1-crossing".
\vskip .2cm 
  Homological degree of the Khovanov complex in the number of 1's in the plane diagram resolution.
  The sum of "quantum" components of the same homological degree i gives the ith component of
  the Khovanov complex.
  
  \vskip .3cm
   One can see that the i-th differential $d^i$ is the sum over "quantum" components, it will map one of the quantum components in
  homological degree i to perhaps several quantum components of homological degree i+1.
  
  \vskip .3cm
  
  Khovanov theory can be considered as a (1+1) dimensional TQFT. The cubes, that are used in
 it's definition come from the TQFT corresponding to the Frobenius algebra defined by $V, m, \Delta$. 
 As we will see later, our constructions will give the interpretation of Khovanov local system as a
 topological D-brane and will suggest  to study the structure of the category of topological D-branes as a {\bf triangulated category}. 

\vskip.2cm

We prove the following important property of  the Khovanov's complex:
\vskip.2cm
{\bf Theorem 1}. Let k denote the kth crossing point of the knot projection $D$, then for any k the Khovanov's complex $C$ decomposes into a sum of two subcomplexes 
$C= C^k_0 \oplus C^k_1$ with matrix differential of the form 

$$ d_C =  \left(\begin{array}{cc}
d_0&d_{0,1}\\
0&d_1 \end{array}\right)$$

\vskip.2cm

{\bf Proof}. Let $C^k_0$ denote the subcomplex of $C$, consisting of vector spaces, which correspond
to the complete resolutions of $D$, having 0 on the kth place. The differential $d_0$ obtained by restricting $d$ only to the arrows between components of $C^k_0$. We define $C^k_1$ the same
way, by restricting to  the complete resolutions of $D$, having 1 on the kth place.

 The only components of the differential, which are not yet used in our decomposition, are the ones
 which change 0-resolution on the kth place of  $C^k_0$ to 1 on the kth place in  $C^k_1$, we denote
 them $d_{0,1}$.

 One can easily see from the definition of the Khovanov's differential (which can be intuitivly described as "all the ways to change 0-resolution in the ith component of the complex to the 1-resolution in the
 (i+1)st component"), that there is no differential mapping ith component of   $C^k_1$ to the   (i+1)st component of $C^k_0$.
\vskip.2cm

\vskip.2cm

\vskip.2cm

{ \sl   Mirror images and adjoints.} 
  Taking the mirror image of the knot
will dualize Khovanov complex. So if we want to invert the cobordism  between two knots,
we should consider the "dual" cobordism between mirror images of these knots.
\vskip .3cm
\newpage
\newpage

 {\sl 3.2.  Reidemeister and Jacobsson moves.}
 
 \vskip .5cm
  A cobordism  (a surface S  embedded into $R^3 \times [0,1]$) between knots $K_0$ and $K_1$ provide a morphism 
  between the corresponding cohomology:
  
  $$F_S :Kh^{i,j}(D_0) \rightarrow Kh^{i,j+\chi(S)}(D_1)$$
  
  where $D_0$ and $D_1$ are diagrams of the knots $K_0$ and $K_1$ and  $\chi(S)$ is the Euler
  characteristic of the surface.
  \vskip .3cm
  We will distinguish between two types of cobordisms - first, corresponding to the wall crossing
  (and changing the type of the knot). And  second, corresponding to nontrivial loops in chambers
which will reflect the dependence of Khovanov homologies on the selfdiffeomorphisms of the
knot, similar to the Reidemeister moves. In this paragraph we will discuss the second type of
cobordisms.
\vskip .3cm

By a surface $S$ in $\R^4$ we mean an oriented, compact
 surface $S,$ possibly with boundary, properly embedded in $\R^3\times 
 [0,1].$ The boundary of $S$ is then a disjoint union
 $$
  \partial S = \partial_0 S \sqcup  - \partial_1 S
$$
 of the intersections of $S$ with two boundary components of 
  $\R^3\times [0,1]$: 
\begin{eqnarray*} 
  \partial_0 S & = & (S\cap \R^3\times \{ 0\}) \\
   - \partial_1 S & = & (S\cap \R^3\times \{ 1\}) 
\end{eqnarray*} 
Note that $\partial_0 S$ and $\partial_1 S$ are oriented links 
 in $\R^3.$ 

The surface $S$ can be represented by a sequence $J$ of plane diagrams 
 of oriented links where every two consecutive diagrams in $J$ are 
 related either by one of the four Reidemeister moves or by one of the four moves 
  {\it birth, death, fusion} described by  Carter-Saito [CS].

   To each Reidemeister 
  move between diagrams $D_0$ and $D_1$ Khovanov [Kh] associates a 
  quasi-isomorphism map of complexes $C(D_0)\rightarrow C(D_1).$

Given a 
representation $J$ of a surface $S$ by a sequence 
 of diagrams, we can associate to $J$  a map 
 of complexes 
$$
  \varphi_J: C(J_0) \to C(J_1) 
$$

Any link cobordism can be described as a one-parameter family 
$D_t, t \in [0,1]$ of planar diagrams, called a {\bf movie}. 
The $D_t$ are link diagrams, except at finitely many singular points which correspond to
 either a Reidemeister move or a Morse modification. 
Away from these points the diagrams for various $t$ are locally isotopic
. Khovanov explained how local moves induce chain 
maps between complexes, hence homomorphisms between 
homology groups. The same is true for planar isotopies. 
Hence, the composition of these chain maps defines a homomorphism 
between the homology groups of the diagrams of  links.

In his paper [Ja] Jacobsson shows that there are knots, s.t.  a {\bf movie} as above
will give a nontrivial morphism of Khovanov homology:

{\bf Theorem [Ja]}  For oriented links $L_0$ and $L_1$, presented by diagrams $D_0$ and $D_1$, 
an oriented link cobordism $\Sigma$ from $L_0$ to $L_1$, defines
a homomorphism $\mathcal{H}(D_0) \rightarrow \mathcal{H}(D_1)$, invariant 
up to multiplication by -1 under ambient isotopy of $\Sigma$ leaving
$\partial \Sigma$ setwise fixed. Moreover, this invariant is non-trivial.

 Jacobsson constructs a family of derived invariants of link cobordisms with the same source 
and target, which are analogous to the classical Lefschetz numbers
of endomorphisms of manifolds.

The Jones polynomial 
appears as the Lefschetz polynomial of the identity cobordism.
\vskip .3cm
 From our perspective the Jacobsson's theorem shows that the Khovanov local system will have nontrivial
 monodromies on the chambers of the space of knots.

\vskip .5cm
{ \sl 3.3. Wall-crossing morphisms.}
\vskip .4cm

 In 3.2 we described what kind of modifications can occur in the cobordism, when we consider
 the "movie" consisting only of manifolds of the same topological type. These modifications implied
 corresponding monodromies of the Khovanov complex.
 
  However, morphisms that are the most important for Vassiliev-type theories are the "wall-crossing"
  morphisms. We will define them now (locally).
  \vskip .3cm
   Consider two complexes $A^\bullet$ and $B^\bullet$ adjacent to the generic wall of the discriminant. Recall, that
   the discriminant is cooriented ( 2.2). If $B^\bullet$ is "right"  via coorientation (or "further in the Ghys metric 
   form the unknot) of $A^\bullet$, then we shift $B^\bullet$'s grading up by one and consider
    $B^\bullet[1]$:
    \vskip 1.5cm
    \begin {center}
    $ A^\bullet | B^\bullet[1]$
    \end{center}
    \vskip 1cm 
    {\bf Note.}  In general, and this will be very important for us in subsequent chapters, if the complex
   $K^\bullet$ is n steps (via the coorientation) away from the unknot, we shift its grading up by n. Thus
   adjacent complexes will have difference in grading  by one (as above), defined by the coorientation.
   \vskip .3cm
    Now we want to understand what happens to the Khovanov complex when we change the kth 
    over-crossing (in the knot  diagram $D$) to an under-crossing.
 
    We will illustrate these changes  on one of the Bar-Natan's trademark diagrams (with his permission)[BN1].

    By "I" we mark the arrow , connecting components of the complex which will exchange places under wall-crossing morphisms when we
    change  over-crossing to under-crossing for the self-intersection point 1. By "II" when we
    do it for point 2 and "III" when we do it for 3:   
    \newpage
    
{\centering\resizebox*{\textwidth}{!}{\parbox{5.85in}{\begin{gather}
  \def\S#1#2{\fbox{
    \!\!\parbox[c]{10mm}{\includegraphics[width=10mm]{figs/#1.eps}}
    \hspace{-3mm}\raisebox{3.5mm}{${\scriptstyle #2}\,$}
    \hspace{-6mm}\raisebox{-3.5mm}{$\scriptscriptstyle #1$} \hspace{-3mm}
  }}
  \def\C#1{\hspace{-2mm}\includegraphics[width=6mm]{figs/C#1.eps}}
  \def\neg{\hspace{-2mm}}
  \begin{array}{c} \xymatrix@R=1.697cm@C=2cm{
    \includegraphics[width=0.7in]{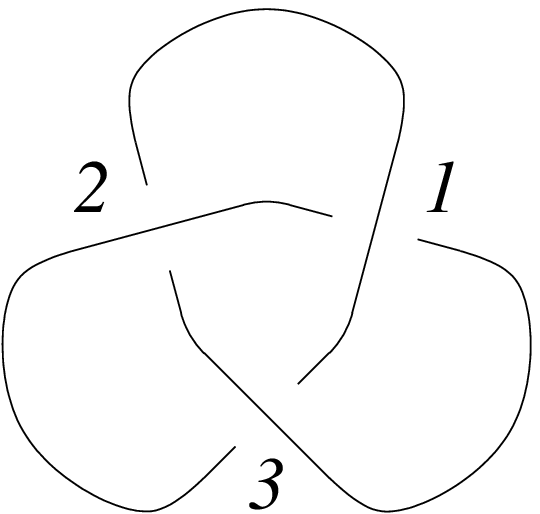}
      & \S{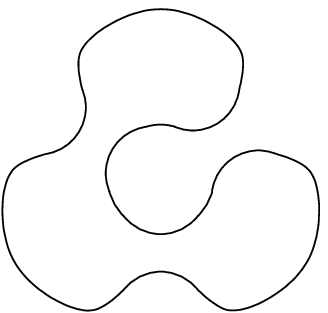}{V\{1\}}
        \ar@{o->}[r]^{\C{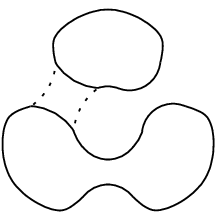}}_{II}
        \ar@{o->}[rd]|\hole^(0.36){\C{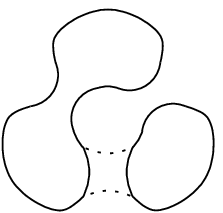}}_(0.6){III}
        \ar@{.}[d]_\oplus
      & \S{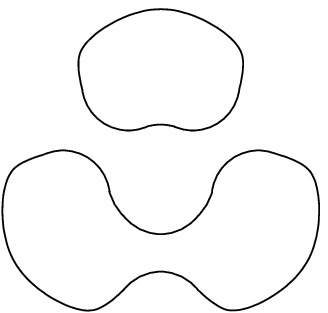}{V^{\otimes 2}\{2\}}
        \ar[rd]^{\C{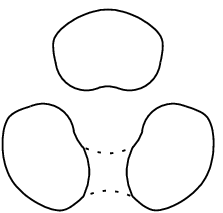}}_{III}
        \ar@{.}[d]^\oplus
      & \\
    \S{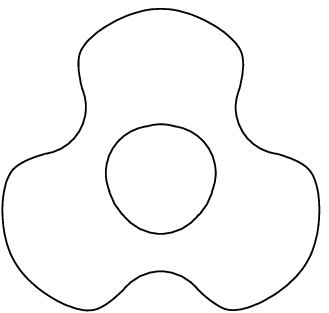}{V^{\otimes 2}}
        \ar[ru]^{\C{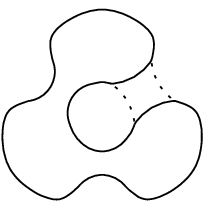}\neg}_{I}
        \ar[r]^{\C{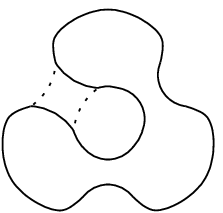}}_{II}
        \ar[rd]^{\C{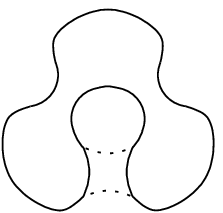}}_{III}
        \ar@{.>}[dd]
      & \S{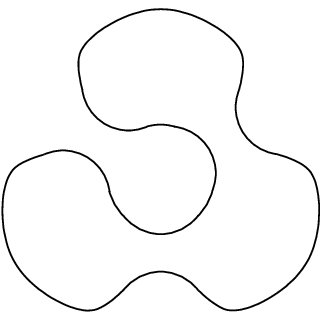}{V\{1\}}
        \ar[ru]^(0.3){\C{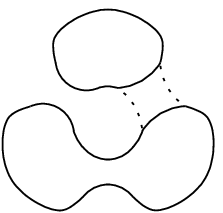}\neg}_(0.66){I}
        \ar@{o->}[rd]^(0.36){\C{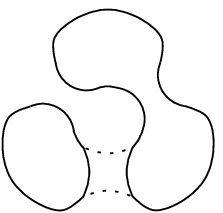}}_(0.6){III}
        \ar@{.}[d]_\oplus
      & \S{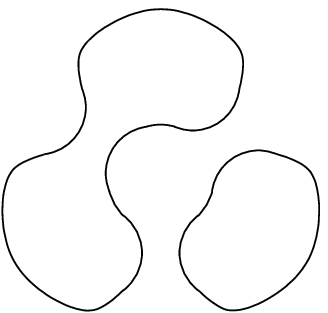}{V^{\otimes 2}\{2\}}
        \ar@{o->}[r]^{\C{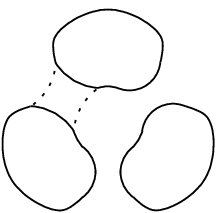}}_{II}
        \ar@{.}[d]^\oplus
      & \S{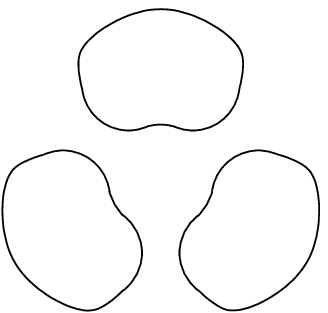}{V^{\otimes 3}\{3\}}
        \ar@{.>}[dd] \\
    \ar@{}[r]_\ ="1"
      & \S{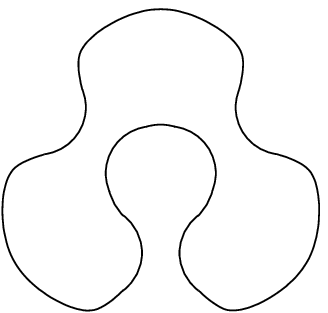}{V\{1\}}
        \ar[ru]|\hole^(0.3){\C{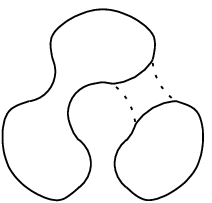}\neg}_(0.66){I}
        \ar[r]^{\C{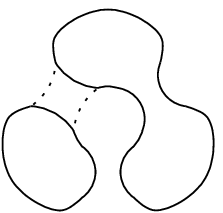}}_{II}="3"
        \ar@{.>}[d]
      & \S{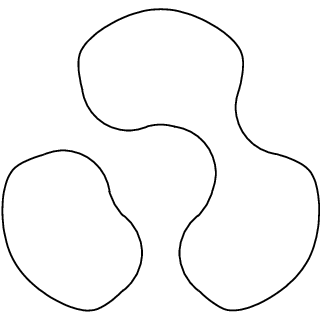}{V^{\otimes 2}\{2\}}
        \ar[ru]^{\C{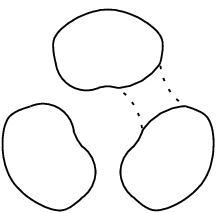}\neg}_{I}
        \ar@{.>}[d]
        \ar@{}[r]_\ ="5" & \\
    \llbracket\righttrefoil\rrbracket^0
        \ar[r]^{d^0}="2"
      & \llbracket\righttrefoil\rrbracket^1
        \ar[r]^{d^1}="4"
      & \llbracket\righttrefoil\rrbracket^2
        \ar[r]^{d^2}="6"
      & \llbracket\righttrefoil\rrbracket^3
  } \end{array}
\notag \\ \label{eq:KCube}
\end{gather}}}}

\vskip .3cm

 Now recall  the theorem proved in (3.1): for any k, where k is the number of  crossings of the diagram $D$, the Khovanov complex can be split into the sum of two subcomplexes with the uppertriangular differential.
 
 Notice from the diagram above that when we change kth overcrossing to an undercrossing, 0 and 1-resolutions are exchanged , so $A^\bullet = A^\bullet_0 \oplus A^\bullet_1$, $B^\bullet[1]=B^\bullet_0[1] \oplus B^\bullet_1[1]$, thus for every k we can define the {\bf wall-crossing morphism}  $\omega$ as follows:
  \vskip 1cm
{\bf Theorem 2.} The map defined as the identity on $ A^\bullet_0 $ and as a trivial map on $A^\bullet_1$:
\begin{equation*}
\xymatrix@C+0.5cm{\omega: A_0^\bullet 
\ar[r]^-{ Id} & B_0^\bullet[1]   \\
\omega:A_1^\bullet \ar[r]^-{ \emptyset} & B_1^\bullet[1] }
\end{equation*}

is the morphism of complexes.
 \vskip 1cm

 {\bf Proof.}  From the Theorem 1 we know that for any crossing k the Khovanov complex can be decomposed as a direct sum with uppertriangular differential:

$$ d =  \left(\begin{array}{cc}
d_0&d_{0,1}\\
0&d_1 \end{array}\right)$$

\vskip.2cm
 It is an easy check that the wall-crossing morphism defined as above is indeed a morphism of complexes (i.e. it commutes with the differential): 
 \vskip 1cm

\begin{equation*}
\xymatrix@C+0.5cm{ A^\bullet
\ar[d]^{d}
\ar[r]^-{ \omega} & B^\bullet \ar[d]^{d}  \\
A^\bullet \ar[r]^-{ \omega} & B^\bullet }
\end{equation*}
 \vskip 1cm

 Since we defined the morphism as 0 on $A_1^\bullet$, the diagram above becomes the following commutative diagram:
\begin{equation*}
\xymatrix@C+0.5cm{ A_0^\bullet
\ar[d]^{d_0}
\ar[r]^-{Id} & B_0^\bullet[1] \ar[d]^{d_0}  \\
A_0^\bullet \ar[r]^-{Id} & B_0^\bullet[1] }
\end{equation*}
 \vskip .5cm

\vskip .3cm
{\sl 3.4. The local system of Khovanov complexes on the space of knots.}
\vskip .5cm

 In this paragraph we introduce the Khovanov local system on the space of knots.
 \vskip .3 cm
 
 {\bf Definition}. A local system on the locally connected topological space M is a fiber bundle over M, the sections of which are abelian groups. The fiber of the bundle depend continuously on the point
 of the base (such that the group structure on the set of fibers can be extended over small domains in the base).
 \vskip .3 cm
 
  Any local system on M with fiber A defines a representation $\pi_1 (M) \rightarrow Aut(A)$. To any
  loop there corresponds a morphism of the fibers of the bundle over the starting point of the loop.
  The set of isomorphism classes of local systems with fiber $A$ are in one-to-one corresondence
  with the set of such representations up to conjugation. For example any representation of an arbitrary
  group $\pi$ in $Aut(A)$ uniquely (up to isomorphism) defines a local system on the space
   $K(\pi, 1)$ [GM].

  Morphisms of local systems are morphisms of fiber bundles, preserving group structure in the fibers.
 Thus introducing the continuation functions (maps between fibers) over paths in the base will define
 a local system over the manifold M.   

  \vskip .3 cm

 Next we set up the Khovanov complexes as a local system on the space of knots. If we were doing it 
 "in coordinates", we would introduce charts on the chambers of the space of knots and define
 our local system via transition maps, starting with some "initial" point . This would be a very interesting and realistic approach, since the homotopy models for chambers are understood [H], [B], e.g. we would have just one chart for the chamber, containing the unknot (since that chamber is contractible), two for a torus knot, four for a hyperbolic one, etc. Then monodromies of the Khovanov local system along nontrivial loops in the chamber will
 be given via Jacobsson movies.

 It would be also very interesting to find a unique special point in every chamber of the space $E$
and study monodromies of the local system with respect to this point. 

The candidate for such point is introduced in the works of J. O'Hara, who studied the minima of the
electrostatical energy  function of the knot [O'H]:

 $$E(K)=\int\int |(x-y)|^{-2} dxdy$$
 
 It was shown  that under some assumptions and for perturbation of the above functional, its critical points  on the space of knots will provide a 'distinguished" point in the chamber. 
 
 The first natural question for this setup is: which nontrivial loops in the
chamber  $E_K$, corresponding to the knot $K$ are distinguised by Khovanov homologies and which are not?
 
 \vskip .5cm
 
 However, assuming  Khovanov's theorem [Kh] (that his homology groups are invariants of the knot,
 independent on the choices made) and assuming also the results of Jacobsson [J] , it is enough for us to introduce the 
 continuation maps, along any path $\gamma$ in the chamber of the space of knots.

  These methods were developed  by several authors (see [Hu]):

 Let $K_1$ and $K_2$ be two knots in the same chamber  of the space $E$, let
${K}_i$ be generic, and let
$\gamma=\{K_t\mid t\in[0,1]\}$ be any path of equivalent objects in
$E$ from $K_1$ to $K_2$.  Then a generic path $\gamma$ 

induces a chain map

$$\label{eqn:continuation}
F({{\gamma}}): CKh_*({K}_1){\longrightarrow}
CKh_*({K}_2)$$

called the ``continuation'' map, which has the following properties:

 \vskip .5cm
\begin{itemize}
\item 
1){\sl Homotopy} A generic homotopy rel endpoints between two
paths ${\gamma}_1$ and ${\gamma}_2$ with
associated chain maps $F_1$ and $F_2$ induces a chain homotopy

$$H:HKh_*({K}_1){\longrightarrow}\nonumber
HKh_{*+1}({K}_2)$$

$$\partial H + H\partial = F_1 - F_2$$
\vskip .3cm
\item
 2){\sl Concatenation} If the final endpoint of
${\gamma}_1$ is the initial endpoint of
${\gamma}_2$, then
$F({{\gamma}_2 {\gamma}_1})$ is chain homotopic
to $F({{\gamma}_2}) F({{\gamma}_1})$.
\vskip .3cm
\item 
3){\sl Constant} If ${\gamma}$ is a constant path
then $F({\gamma})$ is the identity on chains.
\end{itemize}

\vskip .3cm
These three properties imply that if $K_1$ and $K_2$ are equivalent,
then $HKh_*({K}_1)\simeq HKh_*({K}_2)$. (Khovanov's theorem).

 This
isomorphism is generally not canonical, because different homotopy
classes of paths may induce different continuation isomorphisms on
Khovanov homology (Jacobsson moves).  However, since
the loop is contractible, we do know that $HKh_*({K})$
depends only on $K$, so we denote this from now on by $HKh_*(K)$.
\vskip .3cm

   We now define the {\bf restriction} of the Khovanov local system to finite-dimensional subspaces of the 
   space of knots.
   
   Note that in the original setting our complexes may have had different length. For example, the 
   complex
   corresponding to the standard projection of the unknot will have length 1, however, we can
   consider very complicated "twisted" projections of the unknot with an arbitrary large number of 
   crossing points. The corresponding complexes will be quasiisomorphic to the original one.
   
   This construction  resembles the definition of Khovanov homology introduced in [CK], [W]. 
   They define Khovanov 
   homology as a relative theory, where homology groups are calculated
   relative to the twisted unknots.
   
   When considering the restrictions of the Khovanov local system to the subcategories of knots 
   with at most n crossings, we would like {\bf all} complexes to be of length $n+1$.
   
   This can be achieved by "undoing" the local system, starting with the knots of maximal crossing number n
   and then using the wall-crossing morphisms, define complexes of length $(n+1)$, quasiisomorphic to the original ones, in all adjacent chambers. We continue this process till it ends, when we reach the chamber containing unknot.
  
\vskip.3cm 
   
   Recall that Khovanov homology is defined for the knot projection (though is independent of it by Khovanov's theorem).
So we will consider a ramification of Vassiliev space, a pair, the embedding of the circle into $R^3$ and its projection on (x,y)-plane. Then each chamber will be subdivided into "subchambers" corresponding to nonsingular knot projections and the "subdiscriminant" will consist of singular projections of the given knot. The local system, defined on such ramification will live on the universal cover of the base, the original Hatcher chamber corresponding to knot K and morphisms of
the local system between "subchambers" are given by Reidemaister moves. The composition of such moves may constitute the Jacobsson's movie and will give nontrivial monodromies of the local system within the original chamber. 

 {\bf Note.}  As we will see later, if one assines cones of Reidemeister morphisms to the walls of the "subdiscriminant", all such 
cones will be acyclic complexes. This statement in a different form was proved in the original Khovanov [Kh] paper.

\newpage
\centerline{\bf 4. The main definition, invariants of finite type for families.}

\vskip 1cm

{\sl 4.1. Some homological algebra.}

\vskip.5cm
  We describe results and main definitions from the category theory and homological algebra which will be 
used in subsequent chapters. The standard references on this subject are [GM], [Th].

\vskip.2cm
 By constructing the local system of (3.4) we introduced the {\bf derived category} of Khovanov
 complexes. The properties of the derived category are summarized in the axiomatics of the 
 {\bf triangulated category}, which we will discuss in this chapter.
 
 \vskip .3cm
 
{\bf Definition}. An {\sl additive} category is a category $\A$ such that
\begin{itemize}
\item Each set of morphisms $Hom(A,B)$ forms an abelian group.

\item Composition of morphisms distributes over the addition of
  morphisms given by the abelian group structure,
  i.e. $f\circ(g+h)=f\circ g+f\circ h$ and $(f+g)\circ h=f\circ
  h+g\circ h$.

\item There exist products (direct sums) $A\times B$ of any two
  objects $A,B$ satisfying the usual universal properties.

\item There exists a zero object $0$ such that $Hom(0,0)$ is the zero
  group (i.e. just the identity morphism). Thus
  $Hom(0,A)=0=Hom(A,0)$ for all $A$, and the unique zero morphism
  between any two objects is the one that factors through the zero object.
\end{itemize}

\vskip .2cm

So in an abelian category we can talk about exact sequences and
\emph{chain complexes}, and cohomology of complexes. Additive functors
between abelian categories are \emph{exact} (respectively left or
right exact) if they preserve exact sequences (respectively short
exact sequences $0\to A\to B\to C$ or $A\to B\to C\to0$).
\vskip .2cm
{\bf Definition}. The {\sl bounded derived category} $D^b(\A)$ of an abelian category $\A$ has
as objects bounded (i.e. finite length) $\A$-chain complexes, and
morphisms given by chain maps with quasi-isomorphisms inverted as
follows. We introduce morphisms 
$f$ for every chain map between complexes $f:\,X_f\to
Y_f$, and $g^{-1}:\,Y_g\to X_g$ for every quasi-isomorphism
$g:\,X_g\stackrel{\sim\,}{\to}Y_g$. Then form all products of these
morphisms such that
the range of one is the domain of the next. Finally identify any
combination $f_1f_2$ with the composition $f_1\circ f_2$, and
$gg^{-1}$ and $g^{-1}g$ with the relevant identity maps id$_{Y_g}$ and
id$_{X_g}$.
\vskip .2cm

Recall that a triangulated category $C$ is an additive category
equipped with the additional data:
\vskip .2cm
{\bf Definition}. \emph{Triangulated category} is an additive category with a
functor $T: X \rightarrow X [1]$ (where $X^i[1] = X^{i+1}$) and a set of \emph{distinguished
triangles} satisfying a list of axioms. 

The triangles include, for all objects
$X$ of the category:

1) Identity morphism
$$
X \rightarrow X\rightarrow0\rightarrow X[1],
$$

2)
 Any morphism $f:X\rightarrow Y$ can be completed to a distinguished
triangle
$$
X\rightarrow Y\rightarrow C\rightarrow X[1],
$$

3) 
There is also a derived analogue of the 5-lemma, and a compatibility of
triangles known as the octahedral lemma, which can be understood as follows:

If we naively interprete property 1) as the difference $X - X = 0$,
property 2) as $C = X - Y$, then the octahedron lemma says:

$$ (X - Y) - Z =  C - Z = X - (Y - Z)$$

\vskip 1cm

When  topological spaces considered 
up to homotopy there is no notion of kernel or cokernel. The
cylinder construction shows that any map
$f:\,X\to Y$ is homotopic to an inclusion
$X\to\,$cyl$\,(f)=Y\sqcup(X\times[0,1])/f(x)\sim(x,1)$, while the path
space construction shows it is also homotopic to a fibration.

\vskip .5cm

The cone $C_f$ on a map $f:\,X\to Y$ is the space formed from
$Y\sqcup(X\times[0,1])$ by identifying $X\times\{1\}$ with its
image $f(X)\subset Y$, and collapsing $X\times\{0\}$ to a point. 

It  can be considered as a cokernel, i.e. if 
$f:\,X\to Y$

is an inclusion, then $C_f$  is  homotopy
equivalent to $Y/X$.

Taking the $i$th cohomology $H_i$ of each term, and using the
suspension isomorphism $H_i(\Sigma X)\cong H_{i-1}(X)$ gives a sequence

$$H_i(X)\to H_i(Y)\to H_i(Y,X)\to H_{i-1}(X)\to H_{i-1}(Y)\to\ldots$$

which is just the long exact sequence associated to the pair
$X\subset Y$.

Up to homotopy we can make this into a sequence of
simplicial maps, so that taking the associated chain complexes we get
a lifting of the long exact sequence of homology  to the
level of complexes. It exists for all maps $f$, not just
inclusions, with $Y/X$ replaced by $C_f$.

If $f$ is a fibration, $C_f$ can act as the ``kernel'' or
fibre of the map. If $f:\,X\to$\,point, then $C_f=\Sigma X$, the suspension of the fibre $X$.

Thus $C_f$ acts as a combination of both cokernel and kernel,
and if $f:\,X\to Y$ is a map inducing an isomorphism of
homology groups of simply connected spaces then the sequence
$$
H_i(X)\to H_i(Y)\to H_i(C_f)\to H_{i-1}(X)\to H_{i-1}(Y)\to\ldots
$$
implies $H_*(C_f)$=0. Then  $C_f$
homotopy equivalent to a point.  Thus we can give the following definition.
\vskip .3cm

{\bf Definition}. If $X$ and $Y$ are
simplicial complexes, then a simplicial map $f:\,X\to Y$ , defines (up to isomorphism)
an object in triangulated category, called the
{\bf cone of morphism} f, denoted $C_f$.
$$
\renewcommand\arraystretch{1}
C_X^{\bullet}\oplus C_Y^{\bullet}[1] \quad\mathrm{with\ differential}\quad
d_{C_f}=\left(\!\!\!\begin{array}{cc}
d_X & f \\ 0 & d_Y[1] \end{array}\!\!\right),
$$
where $[\,n\,]$ means shift a complex $n$ places up.

Thus we can define the cone $C_f$ on any map of chain complexes
$f:\,A^{\bullet}\to B^{\bullet}$ in an
abelian category $\A$ by the above formula, replacing $C_X^{\bullet}$
by $A^{\bullet}$ and $C_Y^{\bullet}$ by $B^{\bullet}$. If $A^{\bullet}=A$ and $B^{\bullet}=B$
are chain complexes concentrated in degree zero then $C_f$ is the
complex $\{A\stackrel{f\,}{\to}B\}$. This has zeroth cohomology
$h^0(C_f)=$\,ker$\,f$, and $h^1(C_f)=$\,coker\,$f$, so combines the
two (in different degrees). In general it is just the total complex of
$A^{\bullet}\to B^{\bullet}$.

 So what we get in a derived
category is not kernels or cokernels, but ``exact triangles''
$$
A^{\bullet}\to B^{\bullet}\to C^{\bullet}\to A^{\bullet}\,[\,1\,].
$$

Thus we have long exact sequences instead of short exact ones; taking
$i$th cohomology $h^i$ of the above gives the standard long exact sequence
$$
h^i(A^{\bullet})\to h^i(B^{\bullet})\to h^i(C^{\bullet})\to h^{i+1}(A^{\bullet})\to\ldots
$$
 The cone will fit into a triangle:
 
 $$
\xymatrix{
&C\ar[dl]^w_{[1]}&\\
A\ar[rr]^u&&B\ar[ul]^v
} \label{eq:tri1}
$$
The ``$[1]$'' denotes that the map $w$ increases the grade of any object
by one.

\vskip 1cm
{\sl 4.2. Space of knots as a classifying space of the category.}

\vskip .5cm

 In this paragraph we will construct the Khovanov functor from the category of knots into the triangulated
 category of Khovanov complexes.
 \vskip .3cm
{\bf Definition}. The {\bf category of knots} $\mathcal K$  is the category, the objects of which are knots,
$S^1 \rightarrow S^3$, morphisms are cobordisms, i.e. surfaces $ \Sigma$  properly embedded in $\R^3\times [0,1]$ with the boundary of $\Sigma$ 
  being the union of two knots $K_1\subset \R^3\times \{ 0\} $ and 
  $K_2 \subset \R^3\times \{ 1\}$.

  We denote $\mathcal K_n$  the $\bf subcategory$ of knots with at most n crossings.

(Recall that a knot's crossing number is  the lowest number of crossings of any diagram of the knot. )

\vskip .2cm
Note that our cobordisms (morphisms in the category of knots) are {\bf directed} via the coorientation
of the discriminant of the space of knots.

Note that to reverse cobordism, we can consider the same cobordism between mirror images of the knots.
\vskip .2cm

{\bf Definition}.  The {\bf nerve} $\mathcal N (C)$ of a  category C is a simplicial set constructed from the objects and morphisms of C, i.e. points of $\mathcal N (C)$ are objects of $C$, 1-simplices are morphisms of $C$, 2-simplices are commutative triangles, 3-simplices are commutative tetrahedrons of $C$, etc.  $$\mathcal N (C) = (lim  \mathcal N^i (C))$$

 The geometric realization of  a simplicial set   $\mathcal N (C)$ is a topological space, called {\bf the classifying space} of the category C, denoted $B(C)$.
\vskip .3cm

  The following observation is the main guideline for our constructions:
{\sl sheaves  on the classifying space of the category are functors on that category} [Wi].
\vskip .3cm
 Once we prove that the Vassiliev space of knots is a classifying space of the category $\mathcal K$,
 our local system will provide a representation of the Khovanov functor.

 Let C be a category and let Set be the category of sets. For each object A of C let Hom(A,Ð) be the hom functor which maps objects X to the set Hom(A,X).

Recall that a functor $F : C \rightarrow Set$ is said to be {\bf representable} if it is naturally isomorphic to $Hom(A,Ð)$ for some object A of C. A representation of F is a pair $(A, \Psi)$ where

$$   \Psi : Hom(A,Ð) \rightarrow F$$

is a natural isomorphism.

\vskip .5cm
 If $E$ - the space of knots, denote $\mathcal K_E$ the category of knots, whose objects are points in $E$
 and morphisms $Mor (x,y)= \{ \gamma : [0,1] \rightarrow X; s.t. \gamma(0)=x, \gamma(1)=y\}$ and  $\mathcal K_K$ - subcategory corresponding to knots of the same isotopy type K.
\vskip .2cm
{\bf Proposition}. The chamber $E_K$ of the space of long knots for $K$ - unknot, torus of hyperbolic knot is the classifying space of the category $\mathcal K _K$.
\vskip .2cm
{\bf Proof}. By Hatcher's theorem [H] the chambers of the space of knots $E_K$, corresponding to unknot, torus or hyperbolic knot
are $K(\pi, 1)$.


\vskip .2cm

By  definition the space of long knots is $E = \{f: R^1 \rightarrow R^3\}$, nonsingular maps which are standard outside the ball of large radius. 
If $f_1, f_2$ are  vector equations giving knots $K_1, K_2$, then $t f_1 +(1-t)f_2$ is a path in the mapping space, defining a knot for each value of t. The cobordism between two embeddings is given by equations in $R^3\times I$. 
All higher cobordisms can be contracted, since there are no higher homotopy groups in $E_K$. So both the classifying space of the category and the chamber of the space of knots are $K(\pi, 1)$ with the same $\pi$. They are the same as simplicial complexes.

 Note, that in the case of hyperbolic knots one can choose the distinguished point in the chamber - corresponding to the hyperbolic metric on the complement to the knot.

 \vskip .5cm
 
{ \sl 4.3. Vassiliev derivative as a cone of the wall-crossing morphism.}

 \vskip .5cm

    To be able to construct a categorification of Vassiliev theory, we have to extend the local
    system, which we defined on chambers, to the discriminant of the space of knots.
\vskip.2cm

 Recall that according to the axiomatics of the triangulated category, described in (4.1),
  we assign an new object to every morphism in the category:

 for a complex $X=(X^i,d_x^i)$ define a complex 
 
 $X[1]$ by $$(X[1])^i=X^{i+1}, d_{X[1]}=-d_X$$
 
  For a morphism of complexes  $f:X \rightarrow Y$
 let $f[1]:X[1] \rightarrow Y[1]$ coincide with f componentwise.

\vskip.2cm
Let $f:X\rightarrow Y$ be a wall-crossing morphism. The {\bf cone of f} is the following complex $C(f)$:

$$X \rightarrow Y \rightarrow Z=C(f) \rightarrow X[1]$$
 i.e.
$$C(f)^i=X[1]^i \oplus Y^i, d_{C(f)}(x^{i+1},y^i)=(-d_X x^{i+1},f(x^{i+1})-d_Yy^i)$$
\vskip.2cm

 Recall, that we set up the local system on the space of knots (3.4) s.t.  if the complex
   $X^\bullet$ is n steps (via the coorientation) away from the unknot, we shift its grading up by n.  
    So
   complexes in adjacent chambers will have difference in grading  by one, defined by the coorientation.
   
   Thus, given a bigraded complex, associated to the generic wall of the discriminant, we get two natural specialization maps into the neighbourhoods, containing $X^\bullet$ and $Y^\bullet$:
   
  So with any morphism $f$ we associate the triangle:
  
  $$
\xymatrix{
&C_f\ar[dl]^w_{[1]}&\\
X^\bullet \ar[rr]^f&&Y^\bullet \ar[ul]^v
} \label{eq:tri1}
$$

  With any commutative cube
  
  \begin{equation*}
\xymatrix@C+0.5cm{\bullet
\ar[d]^{\omega}
\ar[r]^-{u} & \bullet \ar[d]^{\omega}  \\
\bullet \ar[r]^-{ u} & \bullet }
\end{equation*}
 \vskip 1cm

 (in the space of knots  the above picture corresponds to the  cobordism around the selfintersection of the discriminant of codimension two), we associate the map between cones, corresponding to the
 vertical and horisontal walls, and assign it to the point of their intersection:
  \vskip .5cm

\begin{equation*}
\xymatrix@C+0.5cm{C_u 
\ar[r]^{C_{u \omega}} & C_{\omega}}
\end{equation*}
 \vskip .5cm

 \vskip.2cm
{\bf Lemma}. Given four chambers as above, the order of taking cones of morphisms is irrelevant,
 $C_{u \omega}=C_{\omega u}$.
\vskip.2cm
{\bf Proof}.  see [GM].

    \vskip.3cm
      Consider a point of selfintersection of the discriminant of codimension n. There are $2^n$
   chambers adjacent to this point. Since the discriminant was resolved by Vassiliev [V], this
   point can be considered
    as a point of transversal selfintersection of n hyperplanes in $R^n$, or
   an origin of the coordinate system of $R^n$.
 \vskip.5cm

 Now our local system looks as follows. 
On chambers of our space we have the local system of Khovanov complexes, to  any point $t$ of the generic wall
between chambers containing $X^\bullet$ and $Y^\bullet$  (corresponding to a singular knot),
we assign the cone of the morphism $X^\bullet \rightarrow Y^\bullet$ (with the specialization maps
from the cone to the small neighborhoods of $t$ containing $X^\bullet$ and $Y^\bullet$). To the point
of codimention n we assign the nth cone, $2^n$-graded complex, etc.
\vskip .5cm

{\bf Definition.} The Khovanov homology of the singular knot (with a single double point ) is a bigraded complex 
$$
\renewcommand\arraystretch{1}
X^{\bullet}\oplus Y^{\bullet}[1] \quad\mathrm{with \ the \ matrix \ differential}\quad
d_{C_\omega}=\left(\!\!\!\begin{array}{cc}
d_X & \omega \\ 0 & d_Y[1] \end{array}\!\!\right),
$$
where $X^{\bullet}$ is Khovanov complex of the knot with overcrossing, $Y^{\bullet}$ is the Khovanov complex of the knot with undercrossing and $\omega$ is the wall-crossing morphism. 
\vskip .2cm

 In [S4] we give the geometric interpretation of the above definition.

\vskip.3cm

 {\sl 4.4. The definition of a theory of finite type.} 
  \vskip.5cm

     Once we extended the local system to the singular locus, it is natural to ask if such an extension
  will lead to the categorification of Vassiliev theory.
  
  The first natural guess is that the theory, set up on some space of objects, which has
  quasiisomorphic complexes on all chambers is a theory of order zero. Such theory will
  consist of trivial distinguished triangles as in (a) of the axiomatics of the triangulated category.
 When complexes, corresponding to adjacent chambers are quasiisomorphic, the cone of the morphism
 is an acyclic complex.
 \vskip.3cm

 {\sl Baby example of a theory of order 0.}

 \vskip.2cm
   
  Let M be an n-dimensional compact oriented smooth manifold. Consider the space of functions
  on M. This is an infinite-dimensional Euclidean space. The chambers of the space will
  correspond to Morse functions on M, the walls of the discriminant - to simple degenerations when
  two critical points collide, etc. Let's consider the Morse complex, generated by the critical points
  of a Morse function on M. As it was shown
  by many authors,
  such complex is isomorphic to the CW complex, associated with M. Since we are calculating the homology of M via
  various Morse functions, complexes may vary, but  will have the same homology and Euler characteristics.
  
  Then we can proceed according to our philosophy and assign cones of morphisms to the
  walls and selfintersections of the discriminant. Since complexes on the chambers of the
  space of functions are quasiisomorphic, all cones are acyclic.

  \vskip.2cm
  
   Now we can introduce the main definition of a Floer-type theory  being of finite type n:
   
   \vskip .5cm
  {\bf Main Definition}. The local system of (Floer-type) complexes, extended to the discriminant of the space of manifolds via the cone of morphism, is a {\bf local system of
  order n} if for any selfintersection of the discriminant of codimension $(n+1)$, its (n+1)st cone is an {\bf acyclic complex}.
  \vskip.2cm
  
    How one shows that an 2n-graded complex is acyclic? For example, if one  introduces inverse maps to the
    wall-crossing morphisms
    and construct the homotopy  $\mathcal H$, s.t.:
    
    $$ d \mathcal H - \mathcal H  d = I$$
    
    It is  easy to check  that the existence of such homotopy $\mathcal H$ implies, that the complex doesn't have homology. Suppose $dc=0$, i.e. c is a cycle, then:
    
    $$  d \mathcal H c - \mathcal H  d c =  d \mathcal H c = I$$

    \vskip .5cm
    
    {\bf Example}. Suppose some  local system is conjectured to be of finite type 3. How one would check this? By our definition, we should consider $2^3$ chambers adjacent to the every point of selfintersection
    of the discriminant of codimension 3, and  8 complexes, representing the local system in the small
    neighbourhood of this point. This will correspond to the following commutative cube:

    \vskip 1cm

   \begin{equation*}
\xymatrix@C-0.1cm{ & {B^ \bullet}
\ar[rr]^{h} \ar'[d][dd]^{b} & &
{C^ \bullet} \ar[dd] ^{c}\\
{ A^ \bullet} \ar[ur]_{f} \ar[rr]^(0.65){g}
\ar[dd]^{a} & & {D^ \bullet}
\ar[ur]_{w}
\ar[dd]^{e} & \\
& {F^ \bullet}
\ar'[r]^-{l}[rr] & & {G^ \bullet} \\
{E^ \bullet}
\ar[ru]_{k} \ar[rr]^{m} &&
{H^ \bullet}
\ar[ru]_{n} & }
\end{equation*}
     
 \vskip 1cm 
 Let's write the homotopy equation in the matrix form.

 Consider dual maps   $f^*, g^*,...,w^*$. Then we get formulas for $ d$ and $\mathcal H$  as  $8 \times 8$ matrices:
  \vskip 1cm 
 
 \[
 d =
\left(
\begin{array}{cccccccc}
d_A&f&g&0&a&0&0&0\\
0&d_B&0&h&0&b&0&0\\
0&1&d_D&w&0&0&e&0\\
0&0&1&d_C&0&0&0&c\\
0&0&0&0&d_E&k&m&0\\
0&0&0&0&1&d_F&0&l\\
0&0&0&0&0&1&d_H&n\\
0&0&0&0&0&0&1&d_G\\
\end{array}
\right)
\]

 \vskip 1cm

\[
\mathcal H =
\left(
\begin{array}{cccccccc}
d_A&0&0&1&0&0&0&0\\
f^*&d_B&0&0&0&0&0&0\\
g^*&0&d_D&0&0&0&0&0\\
0&h^*&w^*&d_C&0&0&0&0\\
a^*&0&0&0&d_E&0&0&1\\
0&b^*&0&0&k^*&d_F&0&0\\
0&0&e^*&0&m^*&0&d_H&0\\
0&0&0&c^*&0l^*&h^*&d_G\\
\end{array}
\right)
\]

  \vskip 1cm 
  
   After substituting  these matrices into the equation $d \mathcal H - \mathcal H  d =I$ we obtain  the diagonal matrix which must be homotopic to the identity matrix:
 \vskip 1cm     
 \[
 \left(
\begin{array}{cccccccc}
ff^* + gg^*+ aa^*&0&0&0&0&0&0&0\\
0&-"-&0&0&0&0&0&0\\
0&0&-"-&0&0&0&0&0\\
0&0&0&-"-&0&0&0&0\\
0&0&0&0&-"-&0&0&0\\
0&0&0&0&0&-"-&0&0\\
0&0&0&0&0&0&-"-&0\\
0&0&0&0&0&0&0&cc^* + nn^* + ll^*\\
\end{array}
\right)
\]   
   \vskip 1cm    
  Thus the condition for the local to be of finite type n can be interpreted as  follows. For any selfintersection of the discriminant of codimension $n+1$ consider $2^n$ complexes, forming a commutative cube (representatives of the local system in the chambers
  adjacent to the selfintersection point). Then the naive geometrical interpretation of  the local system 
  being of finite type n is the following: each complex can be
  "split" into $n+1$ subcomplexes, which map quasiisomorphically to $n+1$ neighbours, at
  least no homologies die or being generated.
  
\vskip.2cm
\newpage
{\bf 5.  Knots: theories of finite type. Further directions}.
\vskip .5cm

{\sl 5.1. Examples of combinatorially defined theories.}
 \vskip .5cm
 In the following table we give the examples of theories, which are the categorifications of classical invariants. All these
theories  fit into our framework and may satisfy the finitness condition. 
\vskip .5cm
 
  \begin{center}
     \begin {tabular}{|l|r|}
        \hline
       $\lambda$ & $\lambda = \chi H^* (M)$ \\ \hline
     Jones polynomial &    Khovanov homology [Kh] \\ \hline
    
    Alexander polynomial&  Ozsvath-Szabo knot homology [OS2] \\ \hline
       $sl(n)$ invariants & Khovanov - Rozhansky homology [KhR]\\ \hline
      Casson invariant & Instanton Floer homology [F] \\ \hline
      Turaev's torsion & Ozsvath-Szabo 3 manifold theory [OS1]\\\hline
        Vafa invariant & Gukov-Witten categorification [GW] \\ \hline
      \end{tabular}
  \end{center}

\vskip .7cm

 Note, that the only theory which is not combinatorially defined is the original Instanton Floer homology [F]. The fact that it's Euler characteristics is 
Casson's invariant was proved by C.Taubes [T].

\vskip .5cm
{ \sl 5.2. Generalization to dimension 3 and 4.}
\vskip.5cm
 In our paper [S1] we generalized Vassiliev's construction to the case of 3-manifolds. In [S2] we
 construct the space of parallelizable 4-manifolds and consider the paramentrized version of
 the Refined Seiberg-Witten invariant [BF].

\vskip .2cm
\noindent {\sl a). The space of  3-manifolds and invariants of finite type}.

\vskip .2cm
 Note that all 3-manifolds are parallelizable and therefore
 carry spin-structures.

\vskip .2cm

Following Vassiliev's approach to classification of knots,
  we constructed spaces $E_1$ and  $E_2$  of 3-manifolds by a version of the Pontryagin-Thom construction.

Our main results are as follows:

\vskip .3cm

{\bf Theorem [S1].} In $E_1-D$ each connected component corresponds to a homeomorphism class of 3-dimensional framed manifold.
For any connected framed manifold as above there is one connected component of $E_1-D$
giving its homeomorphism type.

\vskip .2cm

{\bf Theorem [S1].} In $E_2- D$ each connected component
 corresponds to a   homeomorphism class of 3-dimensional spin manifold.
For any connected spin  manifold there is one connected component of
$E_2-D$ giving its homeomorphism type.

\vskip .3cm

By a spin manifold we understand a pair $(M,\theta)$ where  $M$ is an oriented 
3-manifold, and  $\theta$ is a spin structure on $M$. Two spin manifolds $(M,\theta)$ 
and $(M',\theta')$ are called homeomorphic, if there exists a homeomorphism
 $M\rightarrow M'$ taking $\theta$ to $\theta'$. 

\vskip .5cm

  The construction of the space naturally leads to the following definition:

\vskip .4cm
  
{\bf Definition}.  A  map $I: (M,\theta) \rightarrow C$ is called a finite type invariant of (at
 most) order k if it satisfies the condition:
$$\sum_{ L' \in L}(-1)^{\# L'}I(M_{L'})=0$$
where $L'$ is a framed  sublink of link L with even framings, L corresponds to the self-intersection
of the discriminant of codimension k+1, ${\# L'}$ - the number of components of $L'$,  $M_{L'}$ - spin 
3-manifold obtained by surgery on $L'$.

\vskip .3cm

  We introduced an example of Vassiliev invariant of finite order. Given a spin
3-manifold
$M^3$  we consider the Euler characteristic of spin
0-cobordism $W$. Denote  by $I(M,spin) = (sgn (W, spin) -1) (mod2)$.

\vskip .2cm

{\bf Theorem [S1].} Invariant $I(M,spin)$ is finite type of order 1.

\vskip .2cm

The construction of the space of 3-manifolds chambers of which correspond to spin 3-manifolds
is  important for understanding, which additional structures one needs in order to
build the theory of finite-type invariants for homologically nontrivial manifolds. 
It suggests that one should consider spin ramifications of known invariants.

\vskip .5cm

 In the following paper we will generalize our constructions and the main definition to the case of
 3-manifolds. We will construct a local system of Ozsvath-Szabo homologies, extend it to the
 singular locus via the cone of morphism and find examples of theories of finite type.

\vskip .2cm
\noindent {\sl b). Stably parallelizable 4-manifolds.} 
\vskip .2cm

 In this section we modify the previous construction [S1] to get the space of parallelizable 4-manifolds.

 \vskip .2cm
 
By the definition the manifold is parallelizable if it admits the global field of frames, i.e. has
a trivial tangent bundle. In the case of 4-manifolds this condition is equivalent to vanishing of Euler and
the second Stieffel-Whitney class. In particular signature and the Euler characterictic of such manifolds will be 0.

 We will use the theorem of Quinn:
 \vskip .2cm
{\bf Theorem } Any  punctured 4-manifold posesses a smooth structure.
 \vskip .2cm 
  Recall also the result of Vidussi, which states that
manifolds diffeomorphic outside a point have the same Seiberg-Witten 
invariants,so one cannot use them to detect eventual inequivalent smooth 
structures. Thus for the purposes of constructing the family version of the Seiberg-Witten invariants, it will be
sufficient for us to consider ``asymptotically flat'' 4-manifolds, i.e. such that outside the ball $B_R$ of some
large radius R they will be given as the set of common zeros of the system of linear equations (e.g. $f_i(x_1,...x_{n+4})=x_i$
for $i=1,...n$.)

\vskip .3cm
 By Gromov's h-principle any smooth 4-manifold (with all of its smooth structures and metrics) can be obtained as a common set of 
zeros of a system of equations in $R^{\N}$ for sufficiently large N.
\vskip .5cm

{\bf Theorem [S2].} Any parallelizable smooth 4-manifold can be obtained as a set of zeros of n functions
on the trivial (n+4)-bundle over $S^n$. Each manifold will be represented by  $|H^1(M,Z_2) \oplus H^3(M,Z)|$ chambers.
\vskip .5cm

 There is a theory which also 
fits into our template - Ozsvath-Szabo homologies for 3-manifolds, Euler characteristic of which is
Turaev's torsion. It would be interesting to show that this theory is also of finite type or decomposes as
the Khovanov theory.
 
\vskip .2cm
\noindent {\sl c). Ozsvath-Szabo theory as triangulated category}.
\vskip .2cm

 In [S3] we put the theory
developed  by P.Ozsvath and Z. Szabo into the context of homological algebra by considering a
local system of their complexes on the space of 3-manifolds and extending it to the singular locus.
 We show that for the
restricted category the Heegaard Floer complex $CF^{\infty}$ is of finite type one. For other versions
of the theory we will be using the new combinatorial formulas, obtained in [SW].
\vskip .2cm

\vskip .2 cm
  Recall that the categorification is the process of  
  replacing sets with categories, functions with functors, and equations between functions by natural isomorphisms between functors.
 On would hope
 that after establishing this correspondence, homological algebra will provide algebraic structures which one
 should assign to geometrical objects without going into the specifics of a given theory. One can see that this approach
 is very useful in topological category, in particular we will be getting knot and link invariants of Ozsvath and
 Szabo after setting up their local system on the space of 3-manifolds.

\vskip .2cm
{\bf Note} Floer homology can be also considered as  invariants for families, so it would be interesting to connect our work to the one of M.Hutchings [H]. His work can be interpreted as construction of local systems corresponding to various Floer-type
theories on the chambers of our spaces. Then we extend them to the discriminant and classify according
to our definition.

\newpage

\vskip.2cm
{ \sl 6.3. Further directions.}
\vskip.5cm

1. There is a number of immediate questions from the finite-type invariants story:
\vskip.2cm
a). What will substitute the notion of the chord diagram? What is the "basis " in the theories
of finite type?

b). What are the "dimensions" of the spaces of theories of order n?

\vskip.2cm

2. What is the representation-theoretical meaning of the theory of finite type?
\vskip.2cm
a). Is it possible to construct a "universal" knot homology theory in a sense of T.Lee [L] ?

b). Is it possible to rise such a "universal" knot homology theory to the Floer-type theory of
3-manifolds?

\vskip.2cm

3.  There are "categorifications" of other knot invariants: Alexander polynomial [OS], HOMFLY
polynomial [DGR ]. These theories also fit into our setting and it will interesting to show that
they decompose into the series of theories of finite type or that their truncations are of
finite type.

\vskip .2cm
4. The next step in our program [S3] is the construction of the local system of Ozsvath-Szabo homologies
on the space of 3-manifolds introduced in [S1]. We  also plan to raise Khovanov theory to the
homological Floer-type  theory of 3-manifolds.
\vskip .2cm
5. It should be also possible to generalize our program to the study of the diffeomorphism group of a 4-manifold
by considering Gukov-Witten [GW] categorification of Vafa invariant on the moduli space 
constructed in [S2].
\newpage
{\bf 5. Bibliography}.
\vskip .5cm

\vskip .2cm

[BF] Bauer S., Furuta M.,  A stable cohomotopy refinement of Seiberg-Witten invariants: I, II, math.DG/0204340.
\vskip .2cm
[BN1]  Bar-Natan D.,On Khovanov's categorification of the Jones polynomial, math.QA/0201043.
\vskip .2cm
[BN2]  Bar-Natan D.,  Vassiliev and Quantum Invariants of Braids, q-alg/9607001.
\vskip .2cm
[Bu] R. Budney,	 Topology of spaces of knots in dimension 3, math.GT/0506524. 
\vskip .2cm
[CJS] Cohen R., Jones J.,Segal G., Morse theory and classifying spaces, preprint 1995.
\vskip .2cm

[CKV] Champanerkar A., Kofman I., Viro O., Spanning trees and Khovanov homology, preprint.
\vskip .2cm
[D] Donaldson S., The Seiberg-Witten Equations and 4 manifold topology, Bull.AMS, v. 33, 1, 1996.
\vskip.2cm
[DGR] Dunfield N., Gukov S.,  Rasmussen J., The Superpolynomial for Knot Homologies , math.GT/0505662.
\vskip .2cm
[F] Floer A., Morse theory for Lagrangian interesections, J. Differ. Geom. 28 (1988), 513547.
\vskip .2cm
[Fu] Fukaya K., Morse homotopy, $A^{\infty}$-categories and Floer homologies, Proc. of the 1993 
Garc Workshop on Geometry and Topology, v.18 of Lecture Notes
series,p.1-102.Seoul Nat.Univ.,1993.
\vskip .2cm
[G-M] Gelfand S., Manin Yu., Methods of homological algebra, Springer, 1996.
\vskip .2cm
[Gh] Ghys E., Braids and signatures, preprint 2004.
 \vskip .2cm

[Ha] A.Hatcher, Spaces of knots, math.GT/9909095
\vskip .2cm
[HaM] A. Hatcher, D. McCullough,Finiteness of Cassifying Spaces of Relative Diffeomorphism Groups of
 3-manifolds,Geom.Top., 1 (1997)
\vskip .2cm
[Hu] Hutchings M., Floer homology of families 1,  preprint SG/0308115 

\vskip.2cm
[Ja] Jacobsson M., An invariant of link cobordisms from Khovanov homology, Algebraic\&Geometric Topology, v. 4 (2004), 1211-1251.
\vskip .2cm
[Kh] Khovanov M., A Categorification of the Jones Polynomial,  Duke Math. J. 101 (2000), no. 3, 359--426.
\vskip .2cm

[KhR] Khovanov M., Rozansky L.,Matrix factorizations and link homology, math.QA/0401268 

\vskip .2cm
[K] Kontsevich M., Feynmann diagrams and low-dimensional topology. First European Congress of Mathematics, Vol 2(Paris, 1992), Progr.Math.,120, Birkhauser,1994.
\vskip .2cm
[L] Lee T.  An Invariant of Integral Homology 3-Spheres Which Is Universal For All Finite Type Invariants,
q-alg/9601002. 
\vskip .2cm

[MOS]  Manolescu C.,  Ozsvath P., Sarkar S., A combinatorial description of knot Floer homology. , math.GT/0607691.

\vskip .2cm
[O'H] O'Hara J., Energy of Knots and Conformal Geometry , World Scientific Publishing (Jun 15 2000).

\vskip .2cm
[O] Ohtsuki T., Finite Type Invariants of Integral Homology 3-Spheres, J. Knot Theory and its Ramifications 5 (1996).
\vskip .2cm
[OS1] Ozsvath P., Szabo Z., Holomorphic disks and three-manifold invariants: properties and applications, GT/0006194.
\vskip .2cm
[OS2] Ozsvath P., Szabo Z., Holomorphic disks and knot invariants, math.GT/0209056.
\vskip .2cm
[R] D.Ruberman,A polynomial invariant of diffeomorphisms of 4-manifold, geometry and Topology, v.2, Proceeding of the Kirbyfest,
p.473-488, 1999.

\vskip .2cm

[S1] Shirokova N., The Space of 3-manifolds, C.R. Acad.Sci.Paris, t.331, Serie1, p.131-136, 2000.
\vskip .2cm
[S2] Shirokova N., On paralelizable 4-manifolds and invariants for families, preprint  2005.
\vskip .2cm
[S3] Shirokova N., The constructible sheaf of Heegaard Floer Homology on the Space of 3-manifolds, in preparation.
\vskip .2cm
[S4] Shirokova N., The finiteness result for Khovanov homology, preprint 2006. 
\vskip .2cm

[SW] Sarkar S.,  Wang J., A combinatorial description of some Heegaard Floer homologies, math.GT/0607777. 
\vskip .2cm

[T] Taubes C.,Casson's invariant and gauge theory, J.Diff.Geom., 31, (1990), 547-599.
\vskip.2cm

[Th] Thomas, R.P.,Derived categories for the working mathematician, math.AG/0001045  
\vskip .2cm
[V] Vassiliev V., Complements of Discriminants of Smooth Maps, Transl. Math. Monographs 98, Amer. Math.Soc., Providence, 1992.
\vskip .2cm
[Vi] Viro O., Remarks on the definition of Khovanov homology,  math.GT/0202199
\vskip .2cm
[W] Wehrli S.,  A spanning tree model for Khovanov homology, math.GT/0409328.

\vskip .5cm

nadya@math.stanford.edu

\end{document}